# Стохастическая онлайн оптимизация.
# Одноточечные и двухточечные нелинейные многорукие бандиты.
# Выпуклый и сильно выпуклый случаи[1]

*Гасников А.В. (к.ф.-м.н., ИППИ РАН, ПреМоЛаб ФУПМ МФТИ) gasnikov@yandex.ru*

*Крымова Е.А. (к.ф.-м.н., ИППИ РАН) ekkrym@gmail.com*

*Лагуновская А.А. (ИПМ им.М.В. Келдыша РАН, МФТИ) a.lagunovskaya@phystech.edu*

*Усманова И.Н. (ПреМоЛаб ФУПМ МФТИ) ilnura94@gmail.com*

*Федоренко Ф.А. (Кафедра МОУ ФУПМ МФТИ) f.a.fedorenko@gmail.com*

**Аннотация**

В работе предложена безградиентная модификация метода зеркального спуска решения задач выпуклой стохастической онлайн оптимизации. Особенностью постановки является допущение, что реализации значений функции нам доступны с небольшими шумами. Цель данной работы – установить скорость сходимости предложенных методов, и определить, при каком уровне шума, факт его наличия не будет существенно сказываться на скорости сходимости.

**Ключевые слова:** метод зеркального спуска, безградиентные методы, методы с неточным оракулом, стохастическая оптимизация, онлайн оптимизация, многорукие бандиты.

## 1. Введение

Данная работа представляет собой попытку перенесения результатов статьи [1] на онлайн контекст [2] – [10]. А именно, следуя работе [1] рассматривается постановка задачи выпуклой стохастической онлайн оптимизации, в которой на каждом шаге (итерации) вместо градиента можно получить только реализацию значения соответствующего этому шагу функции. При этом допускается, что эта реализация доступна с шумом уровня $\delta$, вообще говоря, не случайной природы. Рассматривается две возможности: на одном шаге (при одной реализации) получать зашумленное значение в одной точке и в двух точках. В первом случае говорят, что рассматривается задача о нелинейных многоруких бандитах (иногда добавляя, одноточечных) [7]. Во втором случае говорят о нелинейных многоруких двухточечных бандитах [7]. Принципиальная разница есть именно при таком переходе [3],

---

[7]. Последующее увеличение числа точек не меняет принципиально картину, соответствующую двум точкам [11].

Основная идея заключается в специальном сглаживании исходной постановки задачи, и использовании метода зеркального спуска [1], [7], [10], [12], [13]. Оригинальной составляющей здесь, в частности, является предложенное в данной статье обобщение этой конструкции на случай наличия шумов. Обратим внимание на условие 1 в раздела 2 (следует сравнить, например, с [1], [14], [15]). Это условие позволило с одной стороны изящно распространить известные оценки на случай, когда есть шумы, см. формулы (2), (3) раздела 2, а с другой стороны это условие хорошо подходит под специфику рассматриваемой в статье постановки (можем получать только зашумленные реализации значений функций), что демонстрируется в раздела 3. Основным результатом работы является теорема 1 раздела 3, в которой результаты статьи [1] переносятся на онлайн контекст.

Во избежание большого количества громоздких выражений, мы опустили часть (наиболее очевидных, но громоздких) выкладок, подробно описав, как они могут быть сделаны. Также в изложении мы не стремились к общности. В частности, для большинства оценок данной работы можно не только выписать точные константы в оценках сходимости в среднем (для этого вполне достаточно написанного в данной статье), но и получить оценки вероятностей больших уклонений. Также можно накладывать более общие требования на классы изучаемых семейств функций, делая константы, характеризующие семейство, не универсальными (одинаковыми для всех шагов), а зависящими от номера шага [3], [16].

Полученные оценки, с учетом известных нижних оценок [2], [7], [9], [17], [18], позволяют говорить о том, что в настоящей работе предложены достаточно эффективные методы, доминирующие в ряде случаев существующие сейчас алгоритмы.

## 2. Метод зеркального спуска для задач стохастической онлайн оптимизации с неточным оракулом

Сформулируем основную задачу стохастической онлайн оптимизации с неточным оракулом. Требуется подобрать последовательность $\{x^k\} \in Q$ ($Q$ – выпуклое множество) так, чтобы минимизировать псевдо регрет [2] – [10]:

(1)  $\text{Regret}_N\left(\{f_k(\cdot)\},\{x^k\}\right) = \frac{1}{N}\sum_{k=1}^{N} f_k(x^k) - \min_{x \in Q} \frac{1}{N}\sum_{k=1}^{N} f_k(x)$

на основе доступной информации

$$\{\nabla_x \tilde{f}_1(x^1, \xi^1); ...; \nabla_x \tilde{f}_{k-1}(x^{k-1}, \xi^{k-1})\}$$

при расчете $x^k$. Причем выполнено **условие**[2]

1. для любых $N \in \mathbb{N}$ ($\Xi^{k-1}$ – сигма алгебра, порожденная $\xi^1$, ..., $\xi^{k-1}$)

$$\sup_{\{x^k = x^k(\xi^1,...,\xi^{k-1})\}_{k=1}^{N}} E\left[\frac{1}{N}\sum_{k=1}^{N}\left\langle E_{\xi^k}\left[\nabla_x f_k(x^k, \xi^k) - \nabla_x \tilde{f}_k(x^k, \xi^k)\big|\Xi^{k-1}\right], x^k - x_*\right\rangle\right] \le \sigma,$$

где $x_*$ – решение задачи

$$\frac{1}{N}\sum_{k=1}^{N} f_k(x) \to \min_{x \in Q},$$

---

[2] В частности, если

$$\left\|\nabla_x \tilde{f}_k(x^k, \xi^k) - \nabla_x f_k(x^k, \xi^k)\right\|_* \le \delta, \quad \max_{x,y \in Q}\|x - y\| \le R,$$

то $\sigma \le \delta R$.

$$E_{\xi^k}\left[\nabla_x f_k(x,\xi^k)\right] = \nabla f_k(x).$$

Здесь случайные величины $\{\xi^k\}$ могут считаться независимыми одинаково распределенными. Онлайновость постановки задачи допускает, что на каждом шаге $k$ функция $f_k(\cdot)$ может выбираться из рассматриваемого класса функций враждебно по отношению к используемому нами методу генерации последовательности $\{x^k\}$. В частности, $f_k(\cdot)$ может зависеть от

$$\{x^1,\xi^1,f_1(\cdot);\ldots;x^{k-1},\xi^{k-1},f_{k-1}(\cdot);x^k\}.$$

Относительно класса функций, из которого выбираются $\{f_k(\cdot)\}$, в данной работе в зависимости от контекста будем предполагать выполненными следующие **условия**:

2. $\{f_k(\cdot)\}$ – выпуклые функции (считаем, что это условие имеет место всегда);
3. $\{f_k(\cdot)\}$ – $\gamma_2$-сильно выпуклые функции в $l_2$;
4. для любых $k=1,\ldots,N$, $x \in Q$

$$E_\xi\left[\left\|\nabla_x \tilde{f}_k(x,\xi)\right\|_*^2\right] \le M^2.$$

Выше (и далее в статье) используется стандартная терминология онлайн оптимизации (см., например, [2], [6], [7], [9]). Однако в отечественной литературе на данный момент имеется определенный дефицит работ по этой (достаточно популярной на западе) тематике. В связи с этим было решено "разбавить" данную статью несколькими простыми примерами, которые позволят лучше прочувствовать смысл используемых понятий.

**Пример 1 (взвешивание экспертных решений, линейные потери).** Рассмотрим задачу взвешивания экспертных решений, следуя [2]. Имеется $n$ различных Экспертов. Каждый Эксперт "играет" на рынке. Игра повторяется $N \gg 1$ раз. Пусть $l_i^k$ – проигрыш (выигрыш со знаком минус) Эксперта $i$ на шаге $k$ ($\left|l_i^k\right| \le M$). На каждом шаге $k$ распределяется один доллар между Экспертами, согласно вектору

$$x^k \in Q = S_n(1) = \left\{x \ge 0: \sum_{i=1}^n x_i = 1\right\}.$$

Потери, которые при этом несем, рассчитываются по потерям экспертов $f_k(x) = \langle l^k, x \rangle$. Целью является таким образом организовать процедуру распределения доллара на каждом шаге, чтобы суммарные потери (за $N$ шагов) были бы минимальны. Допускается, что потери экспертов $l^k$ могут зависеть еще и от текущего хода $x^k$. Установленные далее в этом разделе результаты (формула (2) с $R^2 = \ln n$) позволяют утверждать, что если на каждом шаге можно наблюдать лишь зашумленные проигрыши Экпертов

$$\frac{\partial \tilde{f}_k(x^k,\xi^k)}{\partial x_i} = l_i^k + \xi_i^k + \delta_i^k,$$

где $\{\xi_i^k\}$ – независимые одинаково распределенные случайные величины $\xi_i^k \in N(0,1)$, $\left|\delta_i^k\right| \le \delta$, то существует такой способ действий $x^k\left(\nabla_x \tilde{f}_1(x^1,\xi^1),\ldots,\nabla_x \tilde{f}_{k-1}(x^{k-1},\xi^{k-1})\right)$ (метод зеркального спуска с $\|\ \| = \|\ \|_1$, $d(x) = \ln n + \sum_{i=1}^n x_i \ln x_i$, см. ниже), который позволяет с вероятностью не менее 0.999 после $N$ шагов проиграть лучшему (на этом периоде $k=1,\ldots,N$) Эксперту не более $\mathrm{O}\left((M+1)\sqrt{N \ln n} + \delta N\right)$ долларов, что означает (см. формулу (1))

$$\text{Regret}_N\left(\{f_k(\cdot)\},\{x^k\}\right) = O\left((M+1)\sqrt{\frac{\ln n}{N}} + \delta\right).$$

При $\delta = 0$ эта оценка оптимальна для данного класса задач [2]. □

Опишем метод зеркального спуска для решения задачи (1) (здесь можно следовать огромному числу литературных источников, мы в основном будем следовать работам [16], [19]). Введем норму $\|\ \|$ в прямом пространстве (сопряженную норму будем обозначать $\|\ \|_*$) и прокс-функцию $d(x)$ сильно выпуклую относительно этой нормы, с константой сильной выпуклости $\geq 1$. Выберем точку старта

$$x^1 = \arg\min_{x\in Q} d(x),$$

считаем, что $d(x^1) = 0$, $\nabla d(x^1) = 0$.

Введем брэгмановское "расстояние"

$$V_x(y) = d(y) - d(x) - \langle \nabla d(x), y - x\rangle.$$

Везде в дальнейшем будем считать, что

$$d(x) = V_{x^1}(x) \leq R^2 \text{ для всех } x \in Q.$$

Определим оператор "проектирования" согласно этому расстоянию

$$\text{Mirr}_{x^k}(g) = \arg\min_{y\in Q}\left\{\langle g, y - x^k\rangle + V_{x^k}(y)\right\}.$$

Метод зеркального спуска (МЗС) для задачи (1) будет иметь вид, см., например, [19]

$$x^{k+1} = \text{Mirr}_{x^k}\left(\alpha_k \nabla_x \tilde{f}_k(x^k, \xi^k)\right),\ k = 1, ..., N.$$

Тогда при выполнении условии (2) для любого $u \in Q$, $k = 1, ..., N$ имеет место неравенство, см., например, [19]

$$\alpha_k \langle \nabla_x \tilde{f}_k(x^k, \xi^k), x^k - u\rangle \leq \frac{\alpha_k^2}{2}\|\nabla_x \tilde{f}_k(x^k, \xi^k)\|_*^2 + V_{x^k}(u) - V_{x^{k+1}}(u).$$

Это неравенство несложно получить в случае евклидовой прокс-структуры $d(x) = \|x\|_2^2 / 2$ [20] (в этом случае МЗС для задачи (1) есть просто вариант обычного метода проекции градиента). Разделим сначала выписанное неравенство на $\alpha_k$ и возьмем условное математическое ожидание $E_{\xi^k}\left[\cdot\,|\,\Xi^{k-1}\right]$, затем просуммируем то, что получится по $k = 1, ..., N$, используя условие 1. Затем возьмем от того, что получилось при суммировании, полное математическое ожидание, учитывая условие 4. В итоге, выбирая $u = x_*$, получим при условиях 1, 2, 4, $\alpha_k \equiv \alpha$ [10]

$$N \cdot E\left[\text{Regret}_N\left(\{f_k(\cdot)\},\{x^k\}\right)\right] \leq \frac{V_{x^1}(x_*)}{\alpha} - \frac{E\left[V_{x^{N+1}}(x_*)\right]}{\alpha} + \left(\frac{1}{2}M^2\alpha + \sigma\right)N \leq$$
$$\leq \frac{R^2}{\alpha} + \left(\frac{1}{2}M^2\alpha + \sigma\right)N,$$

выбирая[3]

$$\alpha = \frac{R}{M}\sqrt{\frac{2}{N}},$$

получим

---

[3] Можно получить и адаптивный вариант приводимой далее оценки, для этого потребуется использовать метод двойственных усреднений [10], [19], [20].

(2) $$E\left[\operatorname{Regret}_N\left(\{f_k(\cdot)\},\{x^k\}\right)\right] \le MR\sqrt{\frac{2}{N}} + \sigma;$$

при условиях[4] 1, 3, 4, $\alpha_k \equiv (\gamma_2 k)^{-1}$, $\|\ \| = \|\ \|_2$ [9]

(3) $$E\left[\operatorname{Regret}_N\left(\{f_k(\cdot)\},\{x^k\}\right)\right] \le \frac{M^2}{2\gamma_2 N}(1+\ln N) + \sigma.$$

Оценки (2), (3) являются неулучшаемыми с точностью до мультипликативного числового множителя. Причем верно это и для детерминированных (не стохастических) постановок, в которых нет шумов ($\sigma = 0$), в случае оценки (2) при этом можно ограничиться классом линейных функций [2].

**Пример 2.** Пусть $Q = B_p^n(1)$ – единичный шар в $l_p$ норме. Относительно оптимального выбора нормы и прокс-структуры можно заметить следующее: если $p \ge 2$, то в качестве нормы $\|\ \|$ оптимально выбирать $l_2$ норму и евклидову прокс-структуру. Определим $q$ из $1/p + 1/q = 1$. Пусть $1 \le p \le 2$, тогда $q \ge 2$. Если при этом $q = o(\ln n)$, то оптимально выбирать $l_p$ норму, а прокс-структуру задавать прокс-функцией

$$d(x) = \frac{1}{2(p-1)}\|x\|_p^2.$$

Во всех этих случаях

$$R^2 = \max_{x \in Q} d(x) = \mathrm{O}(1).$$

Для $q \ge \Omega(\ln n)$, выберем $l_a$ норму, где

$$a = \frac{2\ln n}{2\ln n - 1},$$

а прокс-структуру будем задавать прокс-функцией

$$d(x) = \frac{1}{2(a-1)}\|x\|_a^2.$$

В этом случае $R^2 = \mathrm{O}(\ln n)$. Детали см., например, в работах [16], [17]. □

## 3. Одноточечные и многоточечные нелинейные многорукие бандиты

Везде в этом разделе будем считать, что все функции $f_k(x)$ и реализации $f_k(x,\xi)$ определены в $Q_{\mu_0}$ – $\mu_0$-окрестности множества $Q$, и удовлетворяют соответствующим условиям из п. 2 именно в $Q_{\mu_0}$.

Пусть требуется подобрать последовательность $\{x^k\} \in Q$ так, чтобы минимизировать псевдо регрет (1) на основе доступной информации ($m = 1, 2$)

$$\left\{\left\{\tilde{f}_1(x_i^1,\xi^1)\right\}_{i=1}^m; ...; \left\{\tilde{f}_{k-1}(x_i^{k-1},\xi^{k-1})\right\}_{i=1}^m\right\}$$

---
[4] Отметим, что при условии 2, мы еще используем неравенство

$$f_k(x^k) - f(x_*) \le \langle \nabla f_k(x^k), x^k - x_* \rangle$$

при преобразовании левой части неравенства в псевдо регрет, а при условии 3 более точное неравенство

$$2(f_k(x^k) - f(x_*)) \le 2\langle \nabla f_k(x^k), x^k - x_* \rangle - \gamma_2 \|x^k - x_*\|_2^2.$$

при расчете $x^k$. Будем предполагать, что имеет место следующее **условие**

5. для любых $k \in \mathbb{N}$, $i = 1,...,m$, $x_i^k \in Q_{\mu_0}$

$$\left| \tilde{f}_k \left( x_i^k, \xi^k \right) - f_k \left( x_i^k, \xi^k \right) \right| \le \delta,$$

$$E_{\xi^k} \left[ f_k \left( x_i^k, \xi^k \right) \right] = f_k \left( x_i^k \right),$$

$$E_{\xi^k} \left[ \tilde{f}_k \left( x_i^k, \xi^k \right)^2 \right] \le B^2;$$

и, в зависимости от контекста, **условия**

6. для любых $k = 1,...,N$, $x, y \in Q_{\mu_0}$ (далее, как правило, это условие будет использоваться при $r = 2$, исключение сделано в таблице 2)

$$\left| f_k(x, \xi) - f_k(y, \xi) \right| \le M_r(\xi) \|x - y\|_r, \quad M_r = \sqrt{E_\xi \left[ M_r(\xi)^2 \right]} < \infty;$$

7. для любых $k = 1,...,N$, $x, y \in Q_{\mu_0}$

$$\left\| \nabla_x f_k(x, \xi) - \nabla_x f_k(y, \xi) \right\|_2 \le L_2(\xi) \|x - y\|_2, \quad L_2 = \sqrt{E_\xi \left[ L_2(\xi)^2 \right]} < \infty.$$

Введем аналоги $\nabla_x \tilde{f}_k(x, \xi)$ из п. 2 ($\mu \le \mu_0$)

$$\nabla_x \tilde{f}_k(x; e, \xi) := \frac{n}{\mu} \tilde{f}_k(x + \mu e, \xi) e, \text{ (при } m = 1\text{)}$$

$$\nabla_x \tilde{f}_k(x; e, \xi) := \frac{n}{\mu} \left( \tilde{f}_k(x + \mu e, \xi) - \tilde{f}_k(x, \xi) \right) e, \text{ (при } m = 2\text{)},$$

где $e \in RS_2^n(1)$, т.е. случайный вектор $e$ равномерно распределен на сфере радиуса 1 в $l_2$. Считаем, что разыгрывание $e$ происходит независимо ни от чего. Аналогично можно определить незашумленную оценку стохастического градиента $\nabla_x f_k(x; e, \xi)$, убрав в правой части тильды (волны).

Онлайновость постановки задачи допускает, что на каждом шаге $k$ функция $f_k(\cdot)$ может выбираться из рассматриваемого класса функций враждебно по отношению к используемому нами методу генерации последовательности $\{x^k\}$. В частности, $f_k(\cdot)$ может зависеть от

$$\left\{ x^1, \xi^1, f_1(\cdot); ...; x^{k-1}, \xi^{k-1}, f_{k-1}(\cdot) \right\}.$$

Более того, при выборе $f_k(\cdot)$ считается полностью известным наша стратегия. Подчеркнем, что поскольку стратегия рандомизированная, то речь идет об описании этой стратегии, а не о реализации. Это означает, что тому, кто подбирает $f_k(\cdot)$, известно, что $e \in RS_2^n(1)$, но не известно как именно мы его разыграем. Это важная оговорка, если допускать, как и в разделе 2, что на каждом шаге $k$ реализация $e_k \in RS_2^n(1)$ становится известной тому, кто враждебно подбирает $f_k(\cdot)$, то нельзя получить оценку псевдо регрета лучше чем $\Omega(N)$ [3]. Причины этого, связаны с введением рандомизации, и на более простой задаче (линейные одноточечные многорукие бандиты) поясняются, например, в работе [10].

Сгладим исходную постановку с помощью локального усреднения по евклидову шару радиуса $\mu > 0$, который будет выбран позже,

$$f_k^\mu(x,\xi) = E_{\tilde{e}}\left[f_k(x+\mu\tilde{e},\xi)\right],$$
$$f_k^\mu(x) = E_{\tilde{e},\xi}\left[f_k(x+\mu\tilde{e},\xi)\right],$$

где $\tilde{e} \in RB_2^n(1)$, т.е. случайный вектор $e$ равномерно распределен на шаре радиуса 1 в $l_2$. Заменим исходную задачу (1) следующей задачей минимизации

(4) $\quad \text{Regret}_N^\mu\left(\{f_k(\cdot)\},\{x^k\}\right) = \frac{1}{N}\sum_{k=1}^N f_k^\mu(x^k) - \min_{x\in Q}\frac{1}{N}\sum_{k=1}^N f_k^\mu(x).$

Это делается для того, чтобы обеспечить выполнение условия 1 п. 2, см. ниже. Будем считать, что имеют место условия 6, 7 (если условие 7 не выполнено, просто полагаем $L_2 = \infty$). Предположим также, что

$$\min\{M_2\mu, L_2\mu^2/2\} \le \varepsilon/2,$$

т.е.

(5) $\quad \mu \le \max\left\{\dfrac{\varepsilon}{2M_2}, \sqrt{\dfrac{\varepsilon}{L_2}}\right\},$

где $\varepsilon = \varepsilon(N)$ определятся из условия (можно также сказать, что из этого условия определяется $N = N(\varepsilon)$)

$$E\left[\text{Regret}_N\left(\{f_k(\cdot)\},\{x^k\}\right)\right] \le \varepsilon.$$

Из [1] следует, что при условии (5), из

$$\text{Regret}_N^\mu\left(\{f_k(\cdot)\},\{x^k\}\right) \le \varepsilon/2$$

для тех же самых последовательностей $\{f_k(\cdot)\},\{x^k\}$, следует

$$\text{Regret}_N\left(\{f_k(\cdot)\},\{x^k\}\right) \le \varepsilon.$$

Далее сконцентрируемся на минимизации сглаженной версии псевдо регрета (4), контролируя при этом выполнение условия (5).

Введенные выше $\nabla_x \tilde{f}_k(x;e,\xi)$ для задачи (4) удовлетворяют условию 1 с $\sigma$ равным, соответственно,

(6) $\quad \sigma \le E\left[\dfrac{\delta n}{N\mu}\sum_{k=1}^N |\langle e_k, r_k\rangle|\right] \le \dfrac{2\delta R\sqrt{n}}{\mu},\quad$ (при $m=1$)

(7) $\quad \sigma \le E\left[\dfrac{2\delta n}{N\mu}\sum_{k=1}^N |\langle e_k, r_k\rangle|\right] \le \dfrac{4\delta R\sqrt{n}}{\mu},\quad$ (при $m=2$)

где $E[r_k^2] \le 2R^2$, $e_k \in RS_2^n(1)$ – не зависит от $r_k = x^k - x_*$. Оценки (6), (7) следуют из того, что [1], [7], [13]

$$E_e\left[\nabla_x f_k(x;e,\xi)\right] = \nabla f_k^\mu(x,\xi),$$

и из явления концентрации равномерной меры на сфере вокруг экватора (при северном полюсе, заданном вектором $r_k$) [21].

Чтобы можно было воспользоваться оценками (2), (3) раздела 2 осталось в условии 4 раздела 2 оценить константу $M$. Выберем в прямом пространстве норму $l_p$, $1 \le p \le 2$ (см. пример 2 раздела 2). Положим $1/p + 1/q = 1$. При $m=1$ и условии 5 имеем оценки [1]

$$M^2 \le \dfrac{(q-1)n^{1+2/q}B^2}{\mu^2},\quad \text{(при } 2 \le q \le 2\ln n\text{)}$$

$$M^2 \le \frac{4n\ln nB^2}{\mu^2}. \quad (\text{при } 2\ln n < q \le \infty)$$

Наиболее интересны случаи, когда $q=2$, $q=\infty$

(8) $\quad M^2 \le \dfrac{n^2 B^2}{\mu^2}$, (при $q=2$)

(9) $\quad M^2 \le \dfrac{4n\ln nB^2}{\mu^2}$. (при $q=\infty$)

При $m=2$ и выполнении условий условия 5, 6 имеем оценки [1] (случай $2<q<\infty$ рассматривается совершенно аналогично)

$$M^2 \le 3nM_2^2 + \frac{3}{4}n^2 L_2^2\mu^2 + 12\frac{\delta^2 n^2}{\mu^2}, \quad (\text{при } q=2)$$

$$M^2 \le 4\ln nM_2^2 + 3n\ln nL_2^2\mu^2 + 48\frac{\delta^2 n\ln n}{\mu^2}. \quad (\text{при } q=\infty)$$

В частности, если

(10) $\quad \mu \le \min\left\{\max\left\{\dfrac{\varepsilon}{2M_2}, \sqrt{\dfrac{\varepsilon}{L_2}}\right\}, \dfrac{M_2}{L_2}\sqrt{\dfrac{4}{3n}}\right\}, \;\; \delta \le \dfrac{M_2\mu}{\sqrt{12n}}$, (при $q=2$)

(11) $\quad \mu \le \min\left\{\max\left\{\dfrac{\varepsilon}{2M_2}, \sqrt{\dfrac{\varepsilon}{L_2}}\right\}, \dfrac{M_2}{L_2}\sqrt{\dfrac{1}{6n}}\right\}, \;\; \delta \le \dfrac{M_2\mu}{\sqrt{96n}}$, (при $q=\infty$)

то

(12) $\quad M^2 \le 5nM_2^2$, (при $q=2$)

(13) $\quad M^2 \le 5\ln nM_2^2$. (при $q=\infty$)

Далее, полагая в (6), (7), что $\sigma \le \varepsilon/4$, получим дополнительно к (5) (и (10), (11) при $m=2$) условия на $\delta$, $\mu$, $\varepsilon$

$$\frac{2\delta R\sqrt{n}}{\mu} \le \frac{\varepsilon}{4}, \quad (\text{при } m=1)$$

$$\frac{4\delta R\sqrt{n}}{\mu} \le \frac{\varepsilon}{4}, \quad (\text{при } m=2)$$

т.е.

(14) $\quad \delta \le \dfrac{\varepsilon\mu}{8R\sqrt{n}}$, (при $m=1$)

(15) $\quad \delta \le \dfrac{\varepsilon\mu}{16R\sqrt{n}}$. (при $m=2$)

Далее надо воспользоваться оценками (2), (3), добиваясь, соответственно,

(16) $\quad MR\sqrt{\dfrac{2}{N}} \le \dfrac{\varepsilon}{4}$,

(17) $\quad \dfrac{M^2}{2\gamma_2 N}(1+\ln N) \le \dfrac{\varepsilon}{4}$.

Таким образом, при $m=1$ получаем оценки на $\mu(\varepsilon)$ из (5), на $\delta(\varepsilon)$ из (14) и оценки $\mu(\varepsilon)$, на $N(\varepsilon)$ из (8), (9), (16), (17) и оценки $\mu(\varepsilon)$; при $m=2$ получаем оценки на $\mu(\varepsilon)$ из (10), (11), на $\delta(\varepsilon)$ из (10), (11), (15) и оценки $\mu(\varepsilon)$, на $N(\varepsilon)$ из (12), (13), (16), (17).

Не будем здесь выписывать то, что получается – это довольно тривиально, но достаточно громоздко. Вместо этого, резюмируем полученные в работе результаты в более наглядной форме. Для этого введем $\tilde{O}(\ )$. Будем считать, что $\tilde{O}(\ )$ – с точностью до логарифмического множителя (от $n$ и(или) $N$) совпадает с $O(\ )$.

Напомним (обозначения см. в разделе 2 и примере 2), что
$$x^{k+1} = \operatorname{Mirr}_{x^k}\left(\alpha_k \nabla_x \tilde{f}_k\left(x^k; e^k, \xi^k\right)\right), \ k=1,\ldots,N,$$

где $\left\{e^k\right\}_{k=1}^N$ – независимые одинаково распределенные случайные векторы $e^k \in RS_2^n(1)$,

$$\nabla_x \tilde{f}_k\left(x^k; e^k, \xi^k\right) := \frac{n}{\mu} \tilde{f}_k\left(x^k + \mu e^k, \xi^k\right) e^k, \ \text{(при } m=1\text{)}$$

$$\nabla_x \tilde{f}_k\left(x^k; e^k, \xi^k\right) := \frac{n}{\mu}\left(\tilde{f}_k\left(x^k + \mu e^k, \xi^k\right) - \tilde{f}_k\left(x^k, \xi^k\right)\right) e^k, \ \text{(при } m=2\text{)}$$

$$\alpha_k \equiv \alpha = \frac{R}{M}\sqrt{\frac{2}{N}}$$

– в общем случае и

$$\alpha_k \equiv (\gamma_2 k)^{-1},$$

если $f_k(x)$ – $\gamma_2$-сильно выпуклые функции в $l_2$ (в этом случае выбирают $p=2$).

**Теорема 1.** *Пусть рассматривается задача стохастической онлайн оптимизации (1), в постановке, описанной в этом разделе (в безградиентном варианте). Пусть выбрана $l_p$-норма, $1 \le p \le 2$, (см. раздел 2). Согласно этой норме задана прокс-функция и расстояние Брэгмана $V_x(y)$. Пусть $R^2 = V_{x^1}(x_*)$, где $x^1$ и $x_*$ определены в разделе 2. Тогда*

$$E\left[\operatorname{Regret}_{N(\varepsilon)}\left(\left\{f_k(\cdot)\right\}, \left\{x^k\right\}\right)\right] \le \varepsilon,$$

*где $N(\varepsilon)$ определяется в таблицах 1, 2.*

*Таблица 1*

| $m=1$ | $f_k(x)$ – выпуклые функции | $f_k(x)$ – $\gamma_2$-сильно выпуклые функции в $l_2$ норме и $p=2$ |
|---|---|---|
| *Выполнены условия 5, 6* | $\tilde{O}\left(\dfrac{B^2 M_2^2 R^2 n^{1+2/q}}{\varepsilon^4}\right)$ | $\tilde{O}\left(\dfrac{B^2 M_2^2 n^2}{\gamma_2 \varepsilon^3}\right)$ |
| *Выполнены условия 5, 7* | $\tilde{O}\left(\dfrac{B^2 L_2 R^2 n^{1+2/q}}{\varepsilon^3}\right)$ | $\tilde{O}\left(\dfrac{B^2 L_2 n^2}{\gamma_2 \varepsilon^2}\right)$ |

*Таблица 2*

| $m=2$ | $f_k(x)$ – выпуклые функции | $f_k(x)$ – $\gamma_2$-сильно выпуклые функции в $l_2$ норме и $p=2$ |
|---|---|---|
| *Выполнены условия 5, 6* | $\tilde{O}\left(\dfrac{M_p^2 R^2 n^2}{\varepsilon^2}\right)$ | $\tilde{O}\left(\dfrac{M_2^2 n^2}{\gamma_2 \varepsilon}\right)$ |
| *Выполнены условия 5, 6, 7* | $\tilde{O}\left(\dfrac{M_2^2 R^2 n^{2/q}}{\varepsilon^2}\right)$ | $\tilde{O}\left(\dfrac{M_2^2 n}{\gamma_2 \varepsilon}\right)$ |

Обе таблицы заполняются исходя из описанной выше техники. Исключением является вторая строчка таблицы 2, ее мы взяли из [1]. Несложно выписать точные формулы вместо $\tilde{O}(\ )$ во всех полях обеих таблиц. Также несложно выписать условие на допустимый уровень шума $\delta$, при котором мультипликативная константа в точной формуле увеличится, скажем, не более чем в два раза.

Оценки в третьей строчке таблицы 2 неулучшаемы [11] (соответствуют нижним оценкам). Оценки во второй строчке таблицы 2 неулучшаемы по $\varepsilon$ [12], [17]. Все сказанное выше касается и стохастических, но не онлайн постановок [12], [17].

Относительно таблицы 1 имеется гипотеза, что приведенные оценки – неулучшаемы по $n$. По $\varepsilon$ оценки могут быть улучшены за счет ухудшения того, как входит $n$ [18].

В заключение рассмотрим пример, демонстрирующий, что полученные в теореме 1 результаты представляются интересными не только в онлайн контексте.

**Пример 3.** Предположим, что "успешность" некоторого человека зависит от того, как он распоряжается своим временем. Имеется $n$ различных родов деятельности. В $k$-день человек распоряжается своим временем согласно вектору $x^k \in S_n(1)$. Этот вектор отражает доли времени, уделенного соответствующим делам. В конце каждого дня человек получает "обратную связь" от "внешнего мира" вида

$$\tilde{f}(x^k + \mu e^k, \xi^k) = f(x^k + \mu e^k) \cdot (1 + \xi^k),$$

где $\{e^k\}$ – независимые одинаково распределенные случайные векторы $e^k \in RS_2^n(1)$, $\mu$ определяется согласно формуле (5), $\{\xi^k\}$ – независимые (между собой и от $\{e^k\} \in RS_2^n(1)$) одинаково распределенные случайные величины $\xi^k \in N(0,1)$, а выпуклая функция $f(x)$, со свойствами

$$|f(x)| \le B, \quad |f(x) - f(y)| \le M_2 \|x - y\|_2,$$

правильно отражает реальное "положение дел", т.е. минимум этой функции соответствует оптимальной для данного человека конфигурации. Задача человека, заключается в том, чтобы (получая каждый день описанную выше обратную связь), так организовать "процесс своего обучения" (на основе получаемой информации), чтобы как можно быстрее достичь такого состояния[5]

$$\overline{x}^N = \frac{1}{N} \sum_{k=1}^{N} x^k \in S_n(1),$$

что с вероятностью 0.999 имеет место неравенство

$$f(\overline{x}^N) - \min_{x \in S_n(1)} f(x) \le \text{Regret}_N(\{f(\cdot)\}, \{x^k\}) \le \varepsilon.$$

Согласно теореме 1 человек может этого достичь за

$$N = O\left(\frac{B^2 M_2^2 n \ln^2 n}{\varepsilon^4}\right)$$

дней. Если есть возможность получать каждый день информацию $\tilde{f}(x^k + \mu e^k, \xi^k)$ и $\tilde{f}(x^k, \xi^k)$, где $\mu$ определяется согласно формуле (11), то тогда можно улучшить оценку до

---

[5] Заметим, что с достаточной точностью и доверительным уровнем (при $N \gg 1$) можно считать, что

$$\overline{x}^N \simeq \frac{1}{N} \sum_{k=1}^{N} (x^k + \mu e^k).$$

$$N = \mathrm{O}\left(\frac{M_2^2 \ln^2 n}{\varepsilon^2}\right)$$

дней – здесь предполагается также гладкость $f(x)$. □

## 4. Заключение

В работе предложены эффективные методы нулевого порядка (также говорят прямые методы или безградиентные методы) для задач выпуклой стохастической онлайн оптимизации. Методы строились на базе обычного зеркального спуска для задач стохастической оптимизации. Вместо стохастического градиента в зеркальный спуск подставлялись специальные дискретные аналоги, аппроксимирующие стохастический градиент. При правильном пересчете размера шага, получаются эффективные методы.

## Литература